\documentclass{amsart}
\usepackage{latexsym,amsxtra,amscd,ifthen}
\usepackage{amsfonts}
\usepackage{verbatim}
\usepackage{amsmath}
\usepackage{amsthm}
\usepackage{amssymb}

\numberwithin{equation}{section}

\theoremstyle{plain}
\newtheorem{theorem}{Theorem}[section]

\newtheorem{prop}[theorem]{Proposition}

\newtheorem{cor}[theorem]{Corollary}

\theoremstyle{definition}

\newcommand{\N}{\ensuremath{\mathbb N}}

\def\la{\langle}
\def\ra{\rangle}

\def\deg{\mbox{deg\,}}

\def\es{\emptyset}
\def\ext{\mbox{Ext\,}}

\def\gr{\mbox{\bf gr\,}}

\def\gn{{ \mbox{\bf gr\,} Q_n }}  
\def\gq{{ \mbox{\bf gr\,} Q}}

\def\sm{\setminus }

\def\tor{\mbox{Tor\,}}

\begin{document}

\title{Algebras associated to pseudo-roots \\
                of noncommutative polynomials \\
                are Koszul
}

\author{Dmitri Piontkovski}

      \address{ Central Institute of Economics and Mathematics\\
                      Nakhimovsky prosp. 47, Moscow 117418,  Russia}
\thanks{Partially
supported by the grant 02-01-00468 of the Russian Basic Research Foundation}

\email{piont@mccme.ru}

\begin{abstract}

Quadratic  algebras associated to pseudo-roots
                of noncommutative polynomials
                 have been introduced by I. Gelfand, Retakh, and Wilson in
                connection with
                studying the decompositions of noncommutative polynomials.
                Later they (with S. Gelfand and Serconek)
                shown that the Hilbert series of these algebras and their quadratic duals
                satisfy the necessary condition for Koszulity.
                It is proved in this note that these algebras are Koszul.
\end{abstract}



\maketitle

\section{Introduction}

In this note we investigate the quadratic algebras $Q_n$ associated with pseudo-roots of noncommutative polynomials
introduced in~\cite{grw}. This is a homogeneous quadratic algebra over a commutative field $k$
(the results do not depend on the choice of it)
with $2^n -1$ generators ($n \ge 2$). It is shown in~\cite{grw,5} how this algebra encodes an information about
factorizations of noncommutative polynomials of degree $n$ and about noncommutative symmetric functions.

Our purpose is to show that $Q_n$ is Koszul.
It was partially done in~\cite{5}. Namely, the Hilbert series
of $Q_n$ and of its quadratic dual algebra $Q_n^!$ have been calculated, and it occurs that these formal power series satisfy the Froberg relation
$Q_n(z) Q_n^!(-z) = 1$, like Hilbert series of dual Koszul algebras. However, this equality do not imply the Koszulity~\cite{pos, roos}.
Moreover, for general quadratic algebra
$A$ it may happens that there is a Koszul
quadratic algebra $S$ such that $A(z) = S(z)$ and $A^!(z) = S^!(z)$ while $A$ is
not Koszul~\cite{pi1}.

Fortunately, for some particular Koszul algebras $S$  the reverse implication  is true.
Assume that $A$ is filtered by a ${\N}$--graded ordered semigroup $\Gamma$
such that the filtration on the  graded component  $A_n$ is induced by the filtration on
the space of generators $V = A_1$.  Then $A$ has also an ${\N}$--filtration induced by the grading of  $\Gamma$.
Then the associated graded algebra $S = \gr A$ has the same Hilbert series as $A$, and
$A$ is Koszul provided that $S$ is (see~\cite{pp}; idea of the proof: the induced filtration on the bar--complex
leads to a spectral sequence $E_1 = \ext_S (k,k) \Longrightarrow \ext_A (k,k)$).
We  use this fact to establish

\begin{theorem}
\label{t-kos}
\label{main}
The algebra $Q_n$ is  Koszul for every $n \ge 2$.
\end{theorem}

Since the Hilbert function of its dual algebra $Q_n^!$ is a polynomial of degree $n$ (see~\cite{5}), we have

\begin{cor}
The algebra  $Q_n$ has global dimension $n$.
\end{cor}


\section{Koszulity}

Recall that a graded connected $k$--algebra $A$ is called Koszul if the trivial module
$k_A$ admits a linear free minimal resolution, that is, every graded vector space
$\tor_i^A (k,k) \simeq \ext^i_A (k,k) $ is concentrated in degree~$i$.
The standard references on Koszulity are~\cite{pri, bac, lo}; a survey could be found in~\cite{pp}.

\subsection{Associated graded algebra $\gn$}

There are several ways to choose the generators of $Q_n$~\cite{grw,5}; let us fix one choice.
Following~\cite{5}, the algebra  $Q_n$  is generated by the variables $r_A$ where $A$ runs
all nonempty subsets $A \subset I_n = \{ 1, \dots, n \}$.
These generators are connected by the relations
$$
     f_{A; i,j} :=  r_A ( r_{A \setminus i  } - r_{A \setminus j  } ) + (r_{A \setminus i  } - r_{A \setminus j  } ) r_{A \setminus \{ i,j\} }
       - r_{A \setminus i  }^2 +  r_{A \setminus j  }^2 =0
$$
for all  $i, j \in A \subset I_n$ (where we simply write $A \setminus i$ instead of $A \setminus \{ i\}$ etc).
These relations are not linearly independent, because, e.~g., $ f_{A; i,j} = - f_{A; j,i} $.

Recall~\cite{5} that $Q_n$ is filtered by the ordered semigroup $\Gamma$ of monomials on the generators
(the order is degree--lexicographical with $r_A < r_B$ iff either $|A| < |B|$ or, if $|A| = |B| = k$
for some $k$, then  $A = \{ a_1, \dots , a_k  \} $ and $B = \{ b_1, \dots , b_k  \} $ with
$a_1 < \dots < a_{k}, b_1 < \dots < b_{k}$ and for some $t < k$ we have $a_i = b_i$ for all $i< t$ but $a_t > b_t$).
The additional ${\N}$--grading on $\Gamma$     is induced by the grading $\deg r_A = |A|$ on generators.
It induces a ${\N}$--filtration  ${\mathcal N}$ on $Q_n$.
By~\cite{5}, the associated graded algebra $\gn$ has the same generators $r_A, A \subset I_n$, and relations
$$
     h_{A; i,j} :=  r_A ( r_{A \setminus i  } - r_{A \setminus j  } ) =0.
$$

Our next purpose is to prove
\begin{theorem}
\label{t-gr-kos}
Algebra $\gn$ is  Koszul for every $n \ge 2$.
\end{theorem}

By the consideration above, we have
\begin{cor}[Theorem~\ref{t-kos}]
\label{c-t1}
Algebra $Q_n$ is  Koszul for every $n \ge 2$.
\end{cor}

\subsection{Groebner basis of $\gn$}

To be familiar with canonical forms of elements in  $\gn$, we first describe
the Groebner basis of its relations. It was essentially done in~\cite{5}, but in another language.
Recall~\cite{ufn} that a set $f = \{ f_i \}$ of elements of a free algebra ${\mathcal F} = k \la  x_1, \dots , x_n \ra $
with Noetherian order ``$<$'' on monomials is called a {\it Groebner basis} of an ideal $I \subset {\mathcal F} $
(= Groebner basis of the factor algebra ${\mathcal F} /I$) iff $f \subset I$ and for every element $a \in I$ its leading
monomial $\hat a$ has a subword equal to a leading monomial $\hat f_i$ of some $f_i$.
A Groebner basis is called {\it reduced},
if every monomial of the decomposition of $ f_i$  has no subword equal to $\hat f_j$ for $j \ne i$, and
the coefficients of $\hat f_i$ in $f_i$ are unit.

For every nonempty $A \subset I_n$ and every $0 \le t \le |A|$, let us define a monomial $S_A^t$ in the free algebra
$F = k \la r_B | \es \ne B \subset I_n \ra$ as follows : if $A = \{ a_1 , \dots , a_m \}$ with $a_1 < \dot a_m$, then
$S_A^t = r_A r_{A \sm \{ a_1\} } \dots r_{A \sm \{ a_1, \dots, a_t \} } $.

\begin{prop}
\label{p-gb}
 The reduced  Groebner basis of $\gn$ consists of all the elements
$g^t_{A,B} = S_A^t r_{A \sm B} - S_A^{t+1}$ for $t \ge 0$, where  $|B| =  t+1, a_1 \notin B \subset A $.
\end{prop}


It follows that both algebras $\gn$ and $Q_n$ has the same linear basis consisting of monomials
({\it normal} monomials w.~r.~t. the Groebner basis above)
$$
    S_{A^1}^{j_1} \dots S_{A^m}^{j_m},
$$
where for every $1\le i < m$ either $A^{i+1} \not\subset A^i$  or  $|A^{i+1}| \ne |A^i| - j_i-1$.
Notice that similar bases has been constructed in~\cite{5}.

\begin{proof}

Let $I$ be an ideal in $F$ generated by the relations $h_{A; i,j}$ of $\gn$.
Let us show
that all the elements $g^t_{A,B} $  lie in this ideal.

First, every element $g^t_{A,B} $ with $a_1 \in B, t \ge 2$ lies in the ideal $I^t$ generated
by $\{ g^j_{A,B} | j < t, A,B \mbox{ as above} \}$, because $g^t_{A,B} = r_A g^{t-1}_{A \sm a_1, B \sm a_1} $.
By induction, assume that  $I^t \subset I$.
Let $b \in B$, let $B' = B \setminus \{ b \}$, and let
$A' = A  \setminus \{ a_1 \}$.
If $a_1 \notin B$,  we have
$$
g^t_{A,B} = r_A g^{t-1}_{A',B'} + g^{t-1}_{A,B'} ( r_{A' \sm B'}  - r_{A \sm B})
            + S^{t-1}_A   h_{ A \sm B' ; b, a_1 }   \in I.
$$

Obviously, the basis $g = \{ g^t_{A,B}  \}$ is reduced.
Since $h_{A; i,j} = g^0_{A, \{ i \} } - g^0_{A, \{ j \} }$, it generates $I$.
By
Diamond lemma (see~\cite{ufn}), it is sufficient to show that every $s$--polynomial
of $g$ is reducible to zero. Since $\widehat {g^t_{A,B}} = S_A^t r_{A \sm B}$, every such
$s$--polynomial has the form
$$
    s_{A,B,C,p,q} =
       g^p_{A,B} S^{q-1}_{A \sm B \sm a_1} r_{A \sm B \sm C} - S^p_{A} g^q_{A \sm B, C} =
         S^p_{A} S^{q+1}_{A \sm B} - S^{p+1}_{A} S^{q-1}_{A \sm B \sm a_1} r_{A \sm B \sm C} ,
$$
where $a_1 \notin B \subset A, a_1 \notin C \subset A \sm B, |B | = p+1, |C| = q+1$, and
$S^{-1}_{*} = 1$. It is clear that every monomial of the form $r_A r_{A \sm B_1} \dots r_{A \sm B_t}$,
where $B_i \subset A, |B_i| = i$,
is reduced to $S^t_A$ by a sequence of reductions w.~r.~t.  the elements of type $g^i_{A \sm \{ a_1, \dots, a_j \} , B_i }$.
Therefore, both monomials in the representation of $s_{A,B,C,p,q}$    are reduced to the same monomial $S^{p+q+1}_{A} $.
%
%
\end{proof}

The following consequence will be used without notification.
\begin{cor}
\label{c-u}
(a) The ideals $r_C \gn, \es \ne C \subset I_n$ have zero pairwise intersections.

(b) Let $G = \{ u_1, \dots, u_s \}$ be a subset of the generators of $\gn$
such that if $B \subset C \in G$, then $B \in G$,  and let
$U$ be a subalgebra (with unit) in  $\gn$ generated by $G$. Then there is a natural isomorphism of left $U$--modules
$$
      \gn \simeq U \oplus U \otimes \bigoplus_{C \not\subset U} r_C \gn .
$$
\end{cor}

\begin{proof}
Obviuosly, a normal form of an element $r_A a \in F$ is $r_A a'$ for some  linear combination $a'$ of normal monomials.
This proves $(a)$.

To prove $(b)$, it is sufficient to notice
that every monomial $a r_C b$ is normal provided that $a \in U$ is normal, $C \notin G$, and $r_C b$ is normal.
\end{proof}

\subsection{Algebra $\gn$ is Koszul}

\begin{proof}[Proof of Theorem~\ref{t-gr-kos}]

Now, we are ready to prove that the algebra $\gn $ is Koszul.
Consider the following complex $K^n$ of right  $\gn $--modules:
$$
     0 \to    K^n_n \to    K^n_{n-1} \to \dots \to K^n_0 \to 0,
$$
where $K^n_0  =  \overline{\gn} = \sum_{r_C} r_C \gn$
is the augmentation ideal   of $\gn$,   and every $K^n_{t+1}, t \ge 0$ is a free module with generators
$S (A:B)$,  where  $A \ne \es, \min A \not\in B \subset A$, and $|B| = t$. The differential
$d : K^n_1 \to K^n_0 $ sends $S (A: \es)$ to $r_A$, while in higher degrees
it is defined on the symbol $S (A:B)$, where  $A = \{ a_1 , \dots, a_m\}$ and
$B = \{ b_1 , \dots, b_{t} \} $ with $a_i < a_{i+1}, b_j < b_{j+1}$,  by the rule
$$
         d ( S (A:B) ) = \sum_{i=1}^{t} (-1)^{t-1} S(A : B) \otimes ( r_{A \sm B}  -  r_{(A \sm a_1)\sm (B \sm b_i)}).
$$
In fact, it is a standard Koszul complex~\cite{pri} for the quadratic algebra $\gn $, but with non-standard generators.

We claim that the complex $K^n$ is acyclic. Then the complex $K^n \oplus k$  (the addition is in the term $K^n_0$) forms a linear free resolution of the trivial
$\gn$--module $k$, hence, it follows that the algebra   $\gn$ is Koszul.

For completeness, let us also define the algebra $Q_1 = \gq_1$ as a free algebra generated by $r_{\{ 1\} }$
and the complex $K^1 : 0 \to S (  \{ 1\}:\es ) \to r_{\{ 1\} } \gq_1 \to 0$. Obviously, it is acyclic.


We will show that $K^n$ is acyclic by the induction on $n \ge 2$. First, notice that $K^n$ is
a direct sum of the complexes $K^n_A$ for $\es \ne A \subset I_n$, where
the differential submodule  $K^n_A$
is spanned by the generators $S (A:B) $, where $B \subset A$,  and $r_A \in K^n_0 = \gn$.
We have to show that every complex $K^n_A$ is acyclic.

If $|A| = k > 0$, then the complex $K^n_A$ is obviously isomorphic to $K^n_{I_k}$.
Consider a subalgebra $P_k \subset \gn $ generated by $\{ r_A | A \subset I_k \}$.
It is canonically isomorphic to $\gq_k$ by~\cite{grw}. Therefore,
the subcomplex
$\widetilde {K^n_{I_k}}  = r_{I_k} P_k \oplus \sum_{B \subset I_k} S (I_k :B) \otimes P_k$
is isomorphic to $K^k_{I_k}$. (It is indeed a subcomplex, because $d (S (I_k :B)) \in \widetilde {K^n_{I_k}} $.)
If $k < n$, it is acyclic by the induction.  Because of the decomposition
$$
     \gn = P_k \oplus  P_k \otimes \bigoplus_{C \not\subset I_k} r_C  \gn  ,
$$
we have   $K^n_{I_k} \simeq \widetilde {K^n_{I_k}}  \otimes ( k + \bigoplus_{C \not\subset  I_k} r_C  \gn  )$.
Therefore, all complexes $K^n_A$   are acyclic
whenever $|A| < n$.

It remains to show that the complex $K = {K^n_{I_n}}$ is acyclic.
Consider the following filtration on the complex $K$:
$$
    0 = J_{0}^n  \subset \dots  \subset J_n^n = K,
$$
where for $i \ge 2$ the differential submodule $J_i^n$ is  spanned by the set $ T_i^n = \{S(I_n: B) | \max B \le i \}$.
We are going to show that all of them are acyclic. Because of the exact triples
$$
      0 \to J_{i}^n \to   J_{i+1}^n \to   J_{i+1}^n /    J_{i}^n \to 0,
$$
it equivalent to say that all quotients $J_{i+1}^n /    J_{i}^n$ are acyclic.
By the induction arguments, we will assume that all the complexes $J_{i}^t$      are acyclic
for $t < n$.

First, the complex $J_n^n / J_{n-1}^n = K / J_{n-1}^n$ is isomorphic to $K_{I_{n-1}}^n$;
the isomorphism is given by the map $S(I_n : B) \to (-1)^{|B|} S(I_{n-1} : B\sm n)$.
By the above, it is acyclic.

Now, consider the quotient complex  $J_{t}^n /    J_{t-1}^n$ for $t \le n-1$.
It is spanned by the set $U_t^n = \{ S(I_n: B) | \max B = t \}$       with the differential
$$
d S(I_n: B ) =  \sum_{i = 1}^k  (-1)^i  S(I_n: B \sm b_i) \otimes (r_{I_n \sm B} - r_{ I_n\sm 1 \sm (B\sm b_i)} )
$$
(where $B = \{ b_1, \dots , b_k , t\}$ with $ b_1 < \dots < b_k < t$).
We have $d U_t^n \subset U_t^n \otimes N_t$, where $N_t$ is the subalgebra in $\gn$
generated by $\{  r_C | t \notin C \}$.
Let us show that the subcomplex $U_t^n \otimes N_t$ is isomorphic to  the complex
$J_{t-1}^{n-1}$. Consider the injective map $U_t^n \otimes N_t \to K^n_{ I_n \sm t } $
given by $\phi : S(I_n: B ) \mapsto S(I_n \sm t : B \sm t  ) $.
The cyclic permutation $\sigma :=   (n \to (n-1) \to \dots \to t \to n)$
maps the generators of  $K^n_{ I_n \sm t }$ to the set $\{ S(I_{n-1} : B) | n \notin B \} $
and maps $N_t$ to $N_n$. Then the image of the composition $\sigma \phi (U_t^n \otimes N_t)$
is exactly  $J_{t-1}^{n-1}$. By the induction, the last complex is acyclic, and so is the complex
$U_t^n \otimes N_t$.
Because
$$
J_{t}^n /    J_{t-1}^n = U_t^n \otimes \gn \simeq U_t^n \otimes N_t \otimes (k \oplus \bigoplus_{C: t \in C } r_C \gn),
$$
the complex $J_{t}^n /    J_{t-1}^n$ is acyclic as well.

Therefore, all the complexes $J_{t}^n$ are acyclic, thus, the algebra $\gn$ is Koszul.
\end{proof}

\section{Acknowledgment}

I am grateful to Mittag--Leffler Institute for hospitality during preparation this note.

 \section{Addition}

 Another proof of Theorem~\ref{main}, due to   S. Serconek and      R. L. Wilson, has been annonced
 in~\cite{qdet}. According to~\cite{re}, this proof is based on another approach
 and seems to be more complificated.

\end{document}